\documentclass[preprint,review,12pt]{elsarticle}
\usepackage[colorlinks]{hyperref}
\usepackage{color}
\usepackage{graphicx}
\usepackage{graphics}
\usepackage{amssymb}
\usepackage{amsmath}
\usepackage{amsthm}
\usepackage{amsfonts}
\usepackage{mathtools}
\usepackage{geometry}
\usepackage{array}
\journal{Physica D: Nonlinear Phenomena}
\begin{document}
\begin{frontmatter}
\title{Unexpected Linearly Stable Orbits in $3$-Dimensional Billiards}
\author{Hassan Attarchi}
\ead{hassan.attarchi@ucr.edu}
\address{Department of Mathematics, University of California Riverside, CA 92521, USA}
\begin{abstract}
In this work, we construct linearly stable periodic orbits in $3$-dimensional domains with boundaries containing focusing components (small pieces of a sphere) where we place these components arbitrarily far apart. It demonstrates that we cannot directly implement the construction methods of planar billiards with focusing boundaries and chaotic properties into $3$-dimensional billiards.
\end{abstract}
\begin{keyword}
Chaotic Billiards \sep Lyapunov exponents \sep linearly stable \sep astigmatism
\end{keyword}
\end{frontmatter}
\section{Introduction}
From the celebrated work of Sinai \cite{Sin70}, we know that dispersing billiards are $K$-systems, and they have strong mixing properties in entire of the phase space. This fact holds for all dispersing billiards of any dimensions \cite{Sin70,SC87}. The work of Bunimovich \cite{Bun79} showed that the focusing components of the boundary may produce the same mixing result and $K$-property. The reason for this effect in billiards with focusing components is the mechanism of defocusing, where the focusing components of the boundary are sufficiently far apart from each other. In \cite{Bun79}, Bunimovich constructed examples of $2$-dimensional billiards with $K$-property such that their boundary components are non-positive constant curvature curves, where the curvature of boundary components is defined with respect to the inward unit normal vectors of the boundary. In other words, he only used arcs of circles and straight lines to build the examples, such as a billiard in a stadium. In \cite{Bun92}, he introduced a class of curves called ``absolutely focusing mirrors" that each incoming infinitesimal beam of parallel rays focuses after the last reflection from these curves. By using absolutely focusing mirrors in the construction of hyperbolic focusing billiards, one can relax the constant curvature condition of focusing pieces (also, see \cite{Bun10}).\par
In \cite{Woj86}, Wojtkowski introduced a large class of $2$-dimensional billiards with focusing components of the boundary, which have nonvanishing Lyapunov exponents. The focusing pieces of the boundary in his examples are called convex scattering curves, and they satisfy $\frac{d^2r}{ds^2}\leq0$, where $r$ is the radius of curvature at each point of the boundary, and $s$ is the arc length parameter of that component. The defocusing mechanism guarantees that if all convex scattering pieces of the boundary have a distance of at least $\frac{2}{\kappa}$, where $\kappa>0$ is the minimum curvature of these convex scattering components, then there is no stable periodic orbit in that billiard.\par
Here is a natural question: \textit{Do the hyperbolic planar billiards with focusing components have higher dimensional counterparts?} Or more precisely, for some given focusing components (e.g., an absolutely focusing mirror), \textit{Is there any minimum distance that must hold between the focusing components of the boundary to construct a chaotic billiard in dimension $\geq3$?}\par
Wojtkowski, in \cite{Woj90}, contributed the negative answer to these questions. Though Wojtkowski's answer seems correct, the example constructed in \cite{Woj90} to support his answer suffering from a mistake. He showed the required distance between the focusing components (semispheres) in his example must be greater than or equal to some positive constant; however, the reverse inequality is correct there. That means, in his example, we cannot put focusing components arbitrarily far apart and still maintain a stable periodic orbit. In the present paper, a family of examples is constructed, where these examples support the negative answer to the presented questions. Here, by using six (small) pieces of a sphere, a linearly stable periodic orbit is constructed in $\mathbf{R}^3$ such that we can arbitrarily increase the distance between these pieces by changing the reflection angle in these components (without any changes in the radius of these pieces). The reason for this phenomenon in higher-dimensional billiards is Astigmatism \cite{Bun00}. In contrast to $2$-dimensional billiards, the studies of focusing mechanisms in high dimensional billiards are much more subtle \cite{Bun88,BR88}, since it happens in infinitely many directions with different focusing strength (i.e. Astigmatism).\par
The structure of this paper is as follows. In section $2$, the evolution of Jacobi Fields (infinitesimal variation) of billiard orbits is studied in dimension $3$. Section $3$ presents the simplest case of examples, where the free paths between two reflections are longer than the usual free paths in a sphere with the same reflection angle. In section $4$, the distance between focusing components arbitrarily increases by increasing the reflection angle inside these focusing pieces.
\section{Orthogonal Jacobi fields}
A mathematical billiard is a dynamical system that corresponds to the uniform motion of a point particle inside a domain with elastic reflections off the boundary. Thus, billiards are geodesic flows on manifolds with boundaries. Thus, it is natural to consider the Jacobi fields (infinitesimal variation) of billiard orbits.\par
\begin{figure}[h]
	\centering
	\includegraphics[width=10cm]{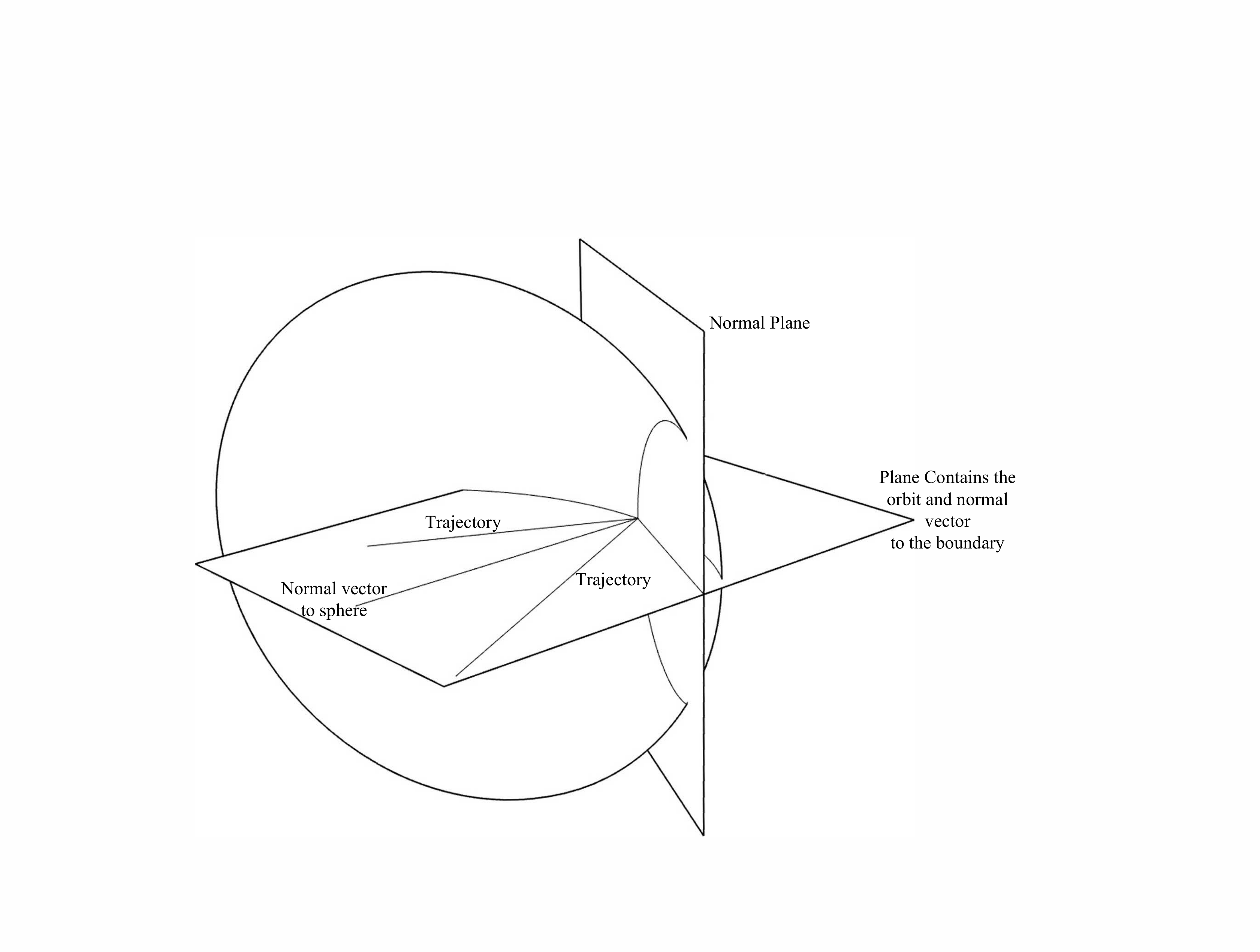}
	\caption{A reflection in $3$-dimensional billiard and the normal plane.}
	\label{0}
\end{figure}

In this work, we consider billiards in domains in $\mathbf{R}^3$, where a Jacobi field of an orbit is a vector field in the plane perpendicular to the orbit at each point (here, we will call this plane the normal plane, see Fig. \ref{0}). Let $J(t)$ be a Jacobi field of an orbit, in this setting $\frac{d^2J}{dt^2}(t)=0$ between reflections when we identify the normal planes by parallel translation. If $L$ denotes the matrix which describes the evolution of Jacobi field $(J(t),\frac{dJ}{dt}(t))$ on the free motion between two consecutive reflections, then
$$L=\begin{bmatrix}
1 & l\\
0 & 1
\end{bmatrix},$$
where $l$ is the length of this free motion.\par
The evolution of a Jacobi field at the moment of the reflection depends on the shape of the boundary. Therefore, it is a complex task to have an explicit formula that describes the evolution of Jacobi fields in all directions in the normal plane at the reflection moment. Here, we will consider billiard orbits inside a domain, where its boundary contains flat components (planes) and pieces of a sphere of radius $r$. Because symmetries exist in these focusing components, it is enough to study the evolution of the Jacobi field in two independent directions.\par
At the reflection moment, the component of the Jacobi field in the plane contains the inbound and outbound trajectory is called the planar component of the Jacobi field. The complement component to the planar component is called the transversal component. If $P$ and $T$ denote the matrices that describe the evolution of the planar and transversal Jacobi fields at reflection moment from a piece of a sphere of radius $r$, respectively, then
\begin{equation}\label{PT}
    P=\begin{bmatrix}
	1 & 0\\
	-\frac{2}{r\cos(\phi)} & 1
	\end{bmatrix},\hspace{1cm}
	T=\begin{bmatrix}
	1 & 0\\
	-\frac{2\cos(\phi)}{r} & 1
	\end{bmatrix},
\end{equation}
where $\phi$ is the reflection angle.
\section{Unexpected linearly stable orbit in 3 dimensional billiards}
\begin{figure}[h]
	\centering
	\includegraphics[width=11cm]{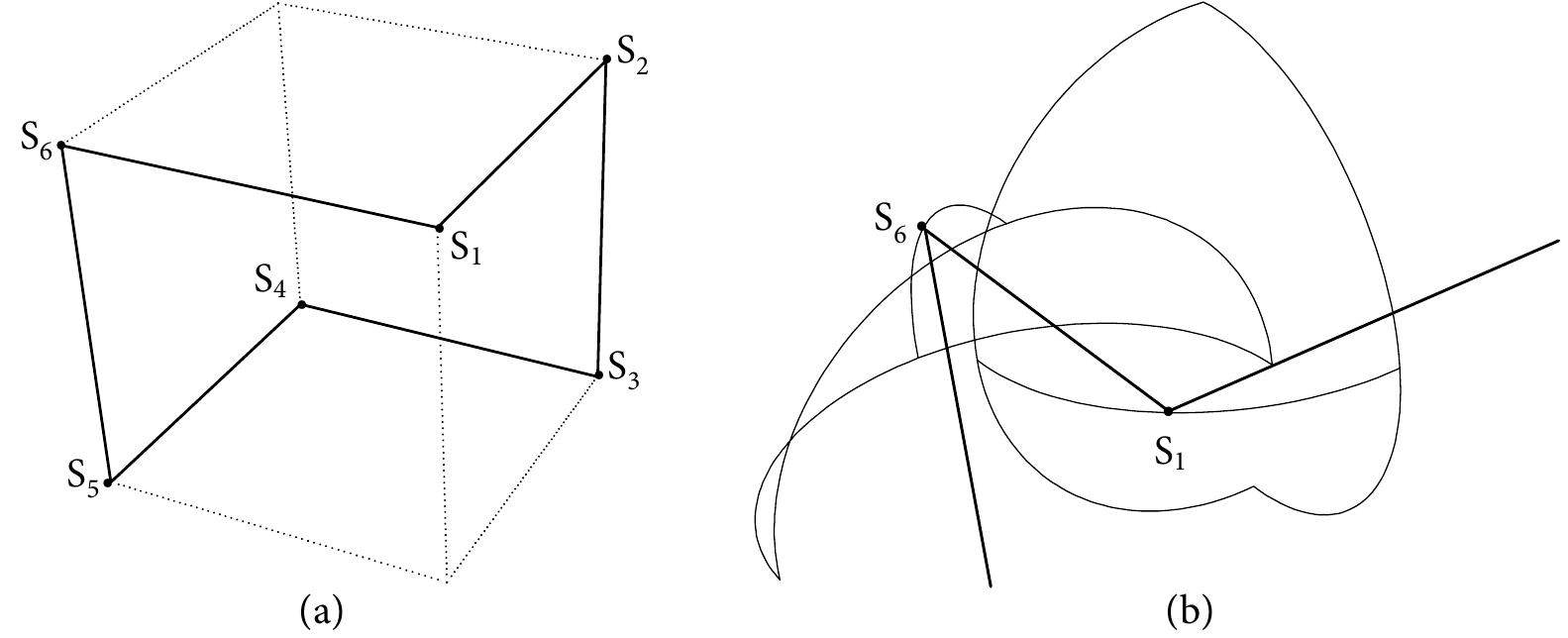}
	\caption{3-dimensional stable orbit: (a) Everything is happening on some edges of a cube. (b) Swapping planar and transversal Jacobi fields in consecutive reflections.}
	\label{1}
\end{figure}

Our periodic orbit has six reflections from some pieces of a sphere of radius $1$, such that the reflection points $\{S_1,S_2,\dots,S_6\}$ are located at some vertices of a cube. And, the periodic orbit stays on some edges of that cube as it is shown in Fig. \ref{1} (a). Therefore, all reflection angles are $\phi=\frac{\pi}{4}$. If we denote the side length of the cube by $l$, then matrices $L$, $P$, and $T$ are respectively given by
\begin{equation}\label{TLP}
L=\begin{bmatrix}
1 & l\\
0 & 1
\end{bmatrix},\hspace{0.5cm} P=\begin{bmatrix}
1 & 0\\
-2\sqrt{2} & 1
\end{bmatrix},\hspace{0.5cm} T=\begin{bmatrix}
1 & 0\\
-\sqrt{2} & 1
\end{bmatrix}.
\end{equation}

In this construction, the planar (transversal) Jacobi field at a reflection becomes the transversal (planar) Jacobi field in the next reflection, see Fig. \ref{1} (b). Hence, the evolution of Jacobi fields around the entire periodic orbit is given by the block diagonal matrix:
$$\begin{bmatrix}
A \ \ & O\\
O\ \  & (TL)A(TL)^{-1}
\end{bmatrix},$$
where
$$A:=(TLPL)^3,$$
and $O$ is a zero block.
Using (\ref{TLP}), we obtain:
$$TLPL=\begin{bmatrix}
1-2\sqrt{2}l & 2(l-\sqrt{2}l^2)\\
-3\sqrt{2}+4l & 1-4\sqrt{2}l+4l^2
\end{bmatrix}.$$
If we consider $A$ as follows
$$A=\begin{bmatrix}
A_{11} & A_{12}\\
A_{21} & A_{22}
\end{bmatrix},$$
then,
\begin{align*}
A_{11} & =1-24\sqrt{2}l+180l^2-224\sqrt{2}l^3+208l^4-32\sqrt{2}l^5,\\
A_{12} & =6l-54\sqrt{2}l^2+272l^3-272\sqrt{2}l^4+224l^5-32\sqrt{2}l^6,\\
A_{21} & =-9\sqrt{2}+156l-360\sqrt{2}l^2+640l^3-240\sqrt{2}l^4+64l^5,\\
A_{22} & =1-30\sqrt{2}l+288l^2-496\sqrt{2}l^3+752l^4-256\sqrt{2}l^5+64l^6,
\end{align*}
and,
\begin{equation}\label{trace}
trace(A)=2-54\sqrt{2}l+468l^2-720\sqrt{2}l^3+960l^4-288\sqrt{2}l^5+64l^6.    
\end{equation}
We are interested in values of $l>0$ such that $-2<trace(A)<2$. Using (\ref{trace}) implies:
$$0<l<\frac{\sqrt{2}}{2}\ \ or\ \ \sqrt{2}<l<\frac{3\sqrt{2}}{2}$$
{\bf Remark:} There are two exception points on each intervals $(0,\frac{\sqrt{2}}{2})$ and $(\sqrt{2},3\frac{\sqrt{2}}{2})$ where $trace(A)=\pm2$.\par
In a sphere with radius $1$ and reflection angle $\phi=\frac{\pi}{4}$ distance between two consecutive reflection points is $\sqrt{2}$. Therefore, by choosing $l\in(\sqrt{2},\frac{3\sqrt{2}}{2})$, we will increase distance between two consecutive reflections, while we maintain the linear stability of the periodic orbit.
\section{Generalization of results}
\begin{figure}[h]
	\centering
	\includegraphics[width=11cm]{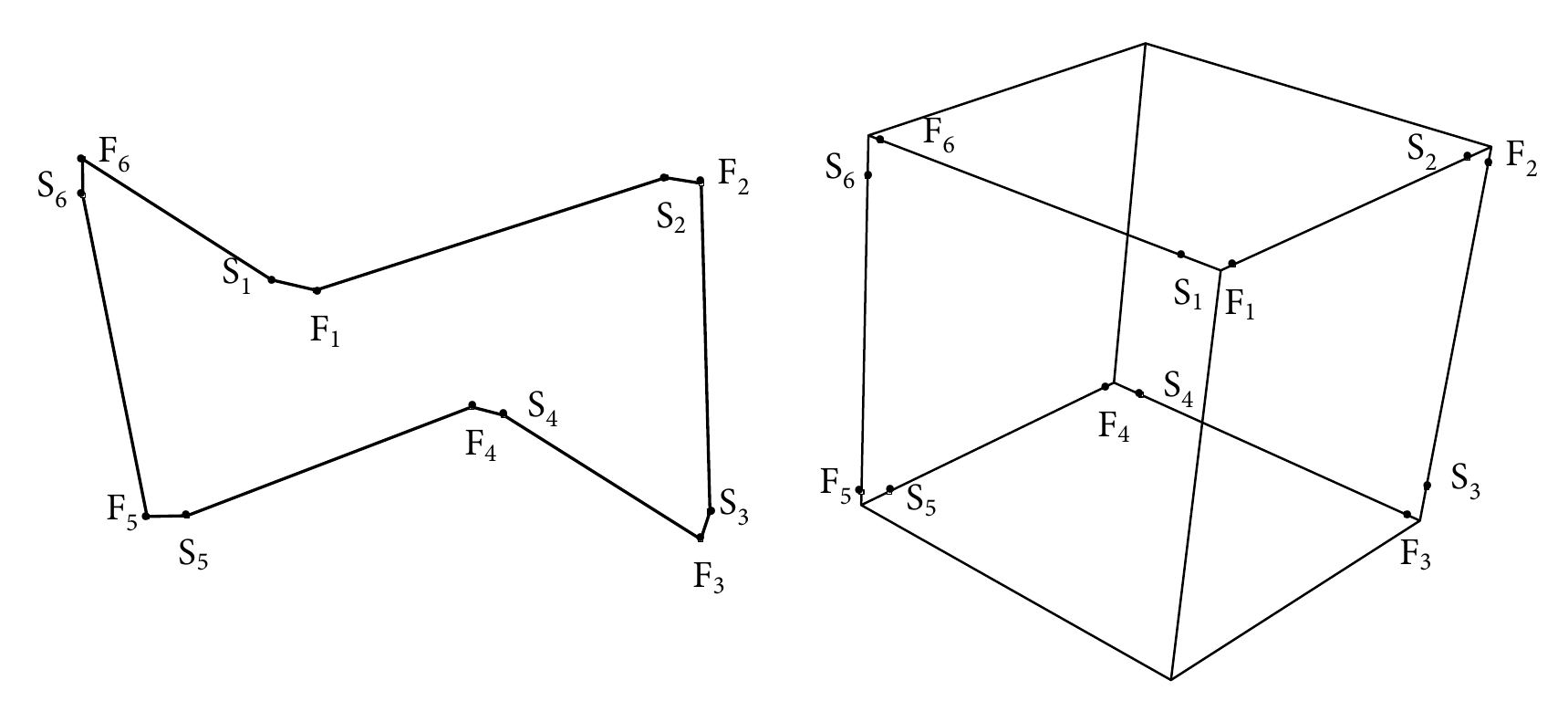}
	\caption{3-dimensional stable orbit with reflections in the following order $S_1,F_1,S_2,F_2,\dots,S_6,F_6$.}
	\label{4}
\end{figure}

In (\ref{PT}), we can see that increasing the reflection angle $\phi$ will decreases the focusing strength in the transversal direction of the Jacobi field at the reflection moment. Thus, one can expect to compensate for this reduction in changes by increasing the distance between reflections. To perform such changes in our variables, we add reflections from some flat components just after reflections from focusing pieces. Therefore, there will be twelve reflections, where six of them are from pieces of a sphere of radius $1$ (they are denoted by $S_1,\dots,S_6$) with reflection angles greater than $\pi/4$. And, there are six other reflections from flat components (they are denoted by $F_1,\dots,F_6$) as it is shown in Fig. \ref{4}. The positions of these flat components are in a way that almost all of the trajectory will remain on some edges of a cube, similar to the example in section 3. Now, we deal with two parameters $l$ (the distance that the point particle travels between two consecutive reflections from focusing components) and $\phi$ (the reflection angle).\par
If $A$ is the matrix that describes the evolution of Jacobi fields around the entire periodic orbit, then
\begin{equation}\label{lphi}
\begin{aligned}
trace(A)= 2 & -36\left(\frac{1}{\cos(\phi)}+\cos(\phi)\right)l+\left(228+96\left(\frac{1}{\cos^2(\phi)}+\cos^2(\phi)\right)\right)l^2\\
& -64\left(\frac{1}{\cos^3(\phi)}+\cos^3(\phi)+\frac{6}{\cos(\phi)}+6\cos(\phi)\right)l^3\\
& +\left(480+192\left(\frac{1}{\cos^2(\phi)}+\cos^2(\phi)\right)\right)l^4\\
& -192\left(\frac{1}{\cos(\phi)}+\cos(\phi)\right)l^5+64l^6.
\end{aligned}
\end{equation}
Let $l:=\frac{1}{\cos(\phi)}$ in (\ref{lphi}), then $trace(A)|_{l=\frac{1}{\cos(\phi)}}=-2$, and
$$\frac{\partial trace(A)}{\partial l}|_{l=\frac{1}{\cos(\phi)}}=36(\frac{1}{\cos(\phi)}-\cos(\phi))>0$$
Thus, for any reflection angles $\pi/4<\phi<\pi/2$, there exists some $\epsilon>0$ such that if $l\in(\frac{1}{\cos(\phi)},\frac{1}{\cos(\phi)}+\epsilon)$, then $-2<trace(A)<2$. That means we can arbitrarily increase the distance between pieces of a sphere in this construction and still have a stable periodic orbit (by increasing the reflection angle $\phi$ towards $\pi/2$).
\section{Conclusion}
The examples of this work show that unlike the billiards with focusing components in $\mathbf{R}^2$, increasing the distance between focusing pieces will not result in instability of periodic orbits in higher dimensional billiards. In dimension $3$, the convergent rate of the Jacobi field varies in different directions in the normal plane after a reflection from a focusing component of the boundary. By mixing the convergent rates of different directions at consecutive reflections, we can balance convergence and divergence to have a stable periodic orbit. Moreover, this work proves that we require some additional conditions to construct chaotic billiards with focusing components in dimension $\geq3$. In \cite{Bun00}, the ``zone of focusing" is defined and used to construct $3$-dimensional billiard tables with focusing components and nonvanishing Lyapunov exponents. It is worth mentioning that the examples and results of this work can easily apply to physical billiards \cite{AB21,Bun19}, where the moving object is a physical particle (a sphere) instead of a mathematical (point) particle.
\section*{Acknowledgement}
I am indebted to L.A. Bunimovich for useful discussions and suggestions.


\begin{thebibliography}{99}
\bibitem{AB21} H. Attarchi and L.A. Bunimovich, Collision of a Hard Ball with Singular Points of the Boundary, Chaos 31, 013123 (2021)
\bibitem{Bun79} L.A. Bunimovich, On the ergodic properties of nowhere dispersing billiards, Commun. Math. Phys. 65, 295-312 (1979)
\bibitem{Bun88} L.A. Bunimovich, Many-dimensional nowhere dispersing billiards with chaotic behavior, Physica D 33, 58-64 (1988)
\bibitem{Bun92} L.A. Bunimovich, On absolutely focusing mirrors, Ergodic Theory and Related Topics III, Springer, Berlin Heidelberg (1992)
\bibitem{BR88} L.A. Bunimovich and J. Rehacek, On the ergodicity of many-dimensional focusing billiards, Annales de l'I.H.P. Physique th\'eorique 68, 421-448 (1998)
\bibitem{Bun00} L.A. Bunimovich, Hyperbolicity and Astigmatism, Journal of Statistical Physics 101, 373-384 (2000)
\bibitem{Bun10} L.A. Bunimovich and A. Grigo, Focusing Components in Typical Chaotic Billiards Should be Absolutely Focusing, Commun. Math. Phys. 293, 127-143 (2010)
\bibitem{Bun19} L.A. Bunimovich, Physical Versus Mathematical Billiards: From Regular Dynamics to Chaos and Back, Chaos 29, 091105 (2019)
\bibitem{Sin70} Ya.G. Sinai, Dynamical systems with elastic reflections, Russ. Math. Surveys 25, 137-189 (1970)
\bibitem{SC87} Ya.G. Sinai and N.I. Chernov, Ergodic properties of some systems of $2$-dimensional discs and
$3$-dimensional spheres, Russ. Math. Surveys 42, 181-207 (1987)
\bibitem{Woj86} M.P. Wojtkowski, Principles for the Design of Billiards with Non-vanishing Lyapunov Exponents, Commun. Math. Phys. 105, 391-414 (1986)
\bibitem{Woj90} M.P. Wojtkowski, Linearly Stable Orbits in $3$ Dimensional Billiards, Commun. Math. Phys. 129, 319-327 (1990)
\end{thebibliography}
\end{document}